\documentclass[10pt]{amsart}

\usepackage{amssymb}
\usepackage{amsmath}
\let\newpf\proof \let\proof\relax 
\newenvironment{pf}{\newpf[\proofname]}{\qed\endtrivlist}

\def\HS{{\mathrm{HS}}}

\def\hull{{\operatorname{hull}}}

\def\bm{\begin{matrix}}
\def\em{\end{matrix}}

\newcommand{\bt}{\begin{thm}}
\newcommand{\et}{\end{thm}}

\newcommand{\bl}{\begin{lemma}}
\newcommand{\el}{\end{lemma}}

\newcommand{\beq}{\begin{eqnarray}}
\newcommand{\eeq}{\end{eqnarray}}

\def\be{\begin{equation}}
\def\ee{\end{equation}}

\def\ba{{\begin{align}}}
\def\ea{{\end{align}}}

\def\SO{{\mathrm {SO}}}

\def\0{{\mathbf 0}}

\def\cal{\mathcal}

\def\SL{{\mathrm {SL}}}

\newtheorem{main thm}{Main Theorem}

\newtheorem{thm}{Theorem}[section]
\newtheorem{cor}[thm]{Corollary}

\newtheorem{lemma}[thm]{Lemma}
\newtheorem{claim}[thm]{Claim}

\theoremstyle{remark}
\newtheorem{rem}{Remark}[section]

\numberwithin{equation}{section}

\def \bn {\hfill \\ \smallskip\noindent}

\theoremstyle{definition}

\def\proof{\bn {\bf Proof.} }

\def\note#1
{\marginpar
{\tiny $\leftarrow$
\par
\hfuzz=20pt \hbadness=9000 \hyphenpenalty=-100 \exhyphenpenalty=-100
\pretolerance=-1 \tolerance=9999 \doublehyphendemerits=-100000
\finalhyphendemerits=-100000 \baselineskip=6pt
#1}\hfuzz=1pt}

\def\tr{{\text{tr}}}

\newcommand{\inter}{\operatorname{int}}
\renewcommand{\mod}{\operatorname{mod}}

\newcommand{\orb}{\operatorname{orb}}

\newcommand{\id}{\operatorname{id}}

\def\loc{{\mathrm{loc}}}

\newcommand{\HH}{{\cal H}}

\newcommand{\C}{{\mathbb C}}

\newcommand{\R}{{\mathbb R}}

\newcommand{\Z}{{\mathbb Z}}

\def\B0{{\bold{0}}}

\catcode`\@=12

\def\Empty{}
\newcommand\oplabel[1]{
  \def\OpArg{#1} \ifx \OpArg\Empty {} \else
  	\label{#1}
  \fi}
		
\newcommand{\comm}[1]{}
\newcommand{\comment}[1]{}

\begin{document}

\title[Limit periodic operators]
{On the spectrum and Lyapunov exponent of limit periodic Schr\"odinger
operators}

\author{Artur Avila}

\date{\today}

\address{
CNRS UMR 7599, Laboratoire de Probabilit\'es et Mod\`eles al\'eatoires\\
Universit\'e Pierre et Marie Curie--Bo\^\i te courrier 188\\
75252--Paris Cedex 05, France
}
\curraddr{IMPA, Estrada Dona Castorina 110, Rio de Janeiro, 22460-320,
Brazil}
\email{artur@math.sunysb.edu}


\begin{abstract}

We exhibit a dense set of limit periodic potentials for
which the corresponding one-dimensional Schr\"odinger operator has
a positive Lyapunov exponent for all energies and a spectrum of zero
Lebesgue measure.  No example with those properties
was previously known, even in the larger class of ergodic potentials.
We also conclude that the generic limit periodic potential has a
spectrum of zero Lebesgue measure.

\end{abstract}

\setcounter{tocdepth}{1}

\maketitle


\section{Introduction}

This work is motivated by a question in the theory of
one-dimensional ergodic Schr\"odinger operators. 
Those are bounded self-adjoint operators of $\ell^2(\Z)$ given by
\be \label {S}
(Hu)_n=u_{n+1}+u_{n-1}+v(f^n(x)) u_n,
\ee
where $f:X \to X$ is an invertible measurable transformation
preserving an ergodic probability measure $\mu$
and $v:X \to \R$ is a bounded measurable function, called the potential.

One is interested in the behavior for
$\mu$-almost every $x$.  In this case, the spectrum is $\mu$-almost surely
independent of $x$.  The Lyapunov exponent is defined as
\be \label {L}
L(E)=\lim \frac {1} {n} \int \ln \|A^{(E)}_n(x)\| d\mu(x),
\ee
where $A^{(E)}_n$ is the
$n$-step transfer matrix of the Schr\"odinger equation $Hu=Eu$.

Here we will give first examples of ergodic potentials with a spectrum of
zero Lebesgue measure such that the Lyapunov exponent is positive
throughout the spectrum.  This answers a question raised by
Barry Simon (Conjecture 8.7 of \cite {S}).

The example we will construct will belong to the
class of limit periodic potentials.  Those arise from continuous potentials
over a minimal translation of a Cantor group (see \S \ref {preliminaries}
for a discussion of those notions).  In our approach, we fix the
underlying dynamics and vary the potential: it turns out that a dense
set of such potentials provide counterexamples.

It is actually possible to incorporate a coupling parameter in our
construction.  Here is a precise version that can be obtained from our
technique:

\begin{thm} \label {coupl}

Let $f:X \to X$ be a minimal translation of a Cantor group.
For a dense set of $v \in C^0(X,\R)$ and for every $\lambda \neq 0$,
the Schr\"odinger operator with potential $\lambda v$ has a spectrum of zero
Lebesgue measure, and the Lyapunov exponent is a continuous positive
function of the energy.

\end{thm}

Our result implies, by continuity of the spectrum, that
a generic potential over a minimal translation of a Cantor group
has a spectrum of zero Lebesgue measure.

\begin{cor} \label {measure}

Let $f:X \to X$ be a minimal translation of a Cantor group.
For generic $v \in C^0(X,\R)$, and for every $\lambda \neq 0$, the
Schr\"odinger operator with potential $\lambda v$ has a spectrum of zero
Lebesgue measure (and the Lyapunov exponent is a continuous function of the
energy which vanishes over the spectrum).

\end{cor}

The statements about the Lyapunov exponent in the generic context are rather
obvious consequences of upper semicontinuity and
density of periodic potentials.  They highlight however that the generic
approach is too rough and that care must be taken
in the proof of Theorem \ref {coupl}
in order not to lose the Lyapunov exponent.

\begin{rem}

Lebesgue measure zero can be strengthened to Hausdorff dimension zero in
both Theorem \ref {coupl} and Corollary \ref {measure}.
It suffices to replace the $10$-th power by an arbitrarily
large one in (2), Lemma \ref {induction}, without
qualitative impact in the proof, and to replace (3), Lemma \ref {joining},
by the covering estimate which the argument is giving.  This remark
was prompted by a recent result (based on a different method)
of Last-Shamis about the Hausdorff dimension
of the spectrum of the critical almost Mathieu operator.

\end{rem}

{\bf Acknowledgements:}  Conjecture 8.7 of \cite {S}
was brought to the attention of the author by Svetlana Jitomirskaya.
This work was carried out
during visits to Caltech and UC Irvine.  This research was partially
conducted during the period the author served as a Clay Research Fellow. 
We are grateful to the referee for several suggestions which led to
significant changes in the presentation.

\section{Preliminaries} \label {preliminaries}

\subsection{From limit periodic sequences to Cantor groups}

Limit periodic potentials are discussed in depth in \cite {AS}.  Here we
will restrict ourselves to some basic facts used in this paper.

Let $\sigma$ be the shift operator on $\ell^\infty(\Z)$,
that is, $(\sigma(x))_n=x_{n+1}$.  Let $\orb(x)=\{\sigma^k(x),\, k \in
\Z\}$.

We say that $x$ is periodic if $\orb(x)$ is finite.  We say that
$x$ is limit periodic if it belongs to the closure, in
$\ell^\infty(\Z)$, of the set of periodic sequences.  If $x$ is limit
periodic, we let $\hull(x)$ be the closure of $\orb(x)$ in
$\ell^\infty(\Z)$.  It is easy to see that every $y \in \hull(x)$ is limit
periodic.

\begin{lemma}

If $x$ is limit periodic then $\hull(x)$ is compact and it has a unique
topological group structure with identity $x$ such that
$\Z \to \hull(x)$, $k \mapsto \sigma^k(x)$ is a homomorphism.
Moreover, the group structure is Abelian and there exist
arbitrarily small compact open neighborhoods of $x$ in
$\hull(x)$ which are finite index subgroups.

\end{lemma}

\begin{pf}

Recall that a metric space is called totally bounded if for every
$\epsilon>0$ it is contained in the $\epsilon$-neighborhood of a finite set.
It is easy to see that a totally bounded subset of a complete
metric space has compact closure.

If $x$ is limit periodic then $\orb(x)$ is
totally bounded: indeed if $p$ is periodic and $\|x-p\|<\epsilon$ then
$\orb(x)$ is contained in the $\epsilon$-neighborhood of $\orb(p)$.
Since $\ell^\infty(\Z)$
is a Banach space, $\hull(x)$ is compact.

Clearly there exists a
unique (cyclic) group structure on $\orb(x)$ such that the
map $\Z \to \orb(x)$, $k \mapsto \sigma^k(x)$ is a homomorphism.

Let us show that the group structure uniformly continuous.
We have
\begin{align} \label {inequalities}
\|\sigma^{k+l}(x)-\sigma^{k'+l'}(x)\|_\infty&=
\|\sigma^{k-k'}(x)-\sigma^{l'-l}(x)\|_\infty\\
\nonumber
&\leq \|\sigma^{k-k'}(x)-x\|_\infty+
\|x-\sigma^{l'-l}(x)\|_\infty\\
\nonumber
&=\|\sigma^k(x)-\sigma^{k'}(x)\|_\infty+
\|\sigma^l(x)-\sigma^{l'}(x)\|_\infty.
\end{align}
where the inequality is just the triangle inequality and the equalities
follow from the fact that $\sigma$ is an isometry of $\ell^\infty(\Z)$.
Thus if
$y,z,y',z' \in \orb(x)$ then $\|y \cdot z-y' \cdot z'\|_\infty \leq
\|y-y'\|_\infty+\|z-z'\|_\infty$, which shows the uniform (even Lipschitz)
continuity.

By uniform continuity, the group structure on $\orb(x)$ has a unique
continuous extension to $\hull(x)$.  Since the group structure on $\orb(x)$
is Abelian, its extension is still Abelian.

For the last statement, fix $\epsilon>0$ and let $p$ be periodic with
$\|x-p\|_\infty<\epsilon/2$.
Let $k$ be such that $\sigma^k(p)=p$.  Clearly
the closure $\hull^k(x)$ of $\{\sigma^{kn}(x),\, n \in \Z\}$ is a compact
subgroup of $\hull(x)$ of index at most $k$.  Since $\hull(x)$ is the union
of finitely many disjoint translates of $\hull^k(x)$, it follows that
$\hull^k(x)$ is also
open.  Since $\sigma$ is an isometry, $\hull^k(x)$ is contained in the
$\epsilon/2$-ball around $p$, and hence it is contained in the
$\epsilon$-ball around $x$.
\end{pf}

By the previous lemma, $\hull(x)$ is compact and totally disconnected, so it
is either finite (if and only if $x$ is periodic) or it is a Cantor set. 

If $x$ is limit periodic but not periodic, we see that every $y$ in
$\hull(x)$ (which is also a limit periodic sequence) is of the form
$y_n=v(f^n(y))$ where $f$ is a minimal translation of a Cantor group
($f=\sigma|\hull(x)$) and $v$ is continuous ($v(w)=w_0$).

\subsection{From Cantor groups to limit periodic sequences}

Let us now consider a Cantor group $X$
and let $t \in X$.  Let $f:X \to X$ be the translation by $t$.  We say
that $f$ is minimal if $\{f^n(y),\, n \in \Z\}$ is dense in
$X$ for every $y \in
X$.  This is equivalent to $\{t^n,\, n \in \Z\}$ being dense in $X$.  In
this case, since there exists a dense cyclic subgroup, we conclude that
$X$ is actually Abelian.

Let $v:X \to \R$ be any continuous function.  Let $\phi:X \to
\ell^\infty(\Z)$, $\phi(x)=(v(f^n(x))_{n \in \Z}$.

\begin{lemma}

For every $x \in X$, $\phi(x)$ is limit periodic and $\phi(X)=\hull(x)$.

\end{lemma}

\begin{pf}

It is enough to show that $\phi(x)$ is limit periodic, since
$\phi(X)$ is compact and $\orb(\phi(x))$ is the image under $\phi$ of the
set $\{f^n(x),\, x \in X\}$ which is dense in $X$.

Given $\delta>0$ we must find a periodic
sequence $p$ such that $\|\phi(x)-p\|_\infty \leq \delta$.
Choose a compact open neighborhood $W$ of the
identity of $X$ which is so small that if $y \in W$ then $|v(y \cdot
z)-v(z)| \leq \delta$.

Introduce a metric $d$ on $X$, compatible with the topology.
Let $\epsilon>0$ be such that if $y,z \in X$ are such that $y \in W$ and $z
\notin W$ then $d(y,z)>\epsilon$.  Choose $m>0$ such that $t^m$ is so close
to the identity that for every $y \in X$, $d(y,f^m(y))<\epsilon$.
Then by induction on $|k|$, $t^{mk} \in W$ for every $k \in \Z$.  It follows
that the closure of $\{t^{km},\, k \in \Z\}$ is a compact subgroup of $X$
contained in $W$.  Clearly it has index at most $m$.

Let $p \in
\ell^\infty(\Z)$ be given by $p_i=v(f^j(w))$ where
$0 \leq j \leq m-1$ is such that $i=j \mod m$.  Then $|\phi(x)_i-p_i|=
|v(f^i(w))-v(f^j(w))|=|v(y \cdot z)-v(z)|$ where $z=f^j(w)$
and $y=t^{i-j}$.  Since $i=j \mod m$, $t^{i-j} \in W$, and by the choice of
$W$ we have $|\phi(x)_i-p_i| \leq \delta$.  It follows that
$\|\phi(x)-p\|_\infty \leq \delta$ as desired.
\end{pf}

\begin{rem}

By the proof above, there exist arbitrarily small
compact subgroups of finite index of $X$ (such subgroups are automatically
open as before).

\end{rem}

\subsection{Limit periodic Schr\"odinger operators}

Given $f:X \to X$ a minimal
translation of a Cantor group and $v:X \to \R$ a continuous function, we
define for every $x \in X$ a
Schr\"odinger operator $H=H_{f,v,x}$ by (\ref {S}).  A formal
solution of $Hu=Eu$ satisfies
\be
A^{(E,f,v)}_n(x) \left (\bm u_0 \\ u_{-1} \em \right )=
\left (\bm u_n \\ u_{n-1} \em \right )
\ee
where
\be
A^{(E)}_n(x)=
A^{(E,f,v)}_n(x)=
S_{n-1} \cdots S_0 \text { where } S_i=\left (\bm E-v(f^i(x))
& -1 \\ 1 & 0 \em \right ).
\ee
The $A^{(E)}_n(x)$ are thus in $\SL(2,\R)$, and are
called the $n$-step transfer matrices.
The Lyapunov exponent $L(E)=L(E,f,v)$
is defined by (\ref {L}), where we take $\mu$ the Haar
probability measure on $X$ (this is the only possible choice actually, since
minimal translations of Cantor groups are uniquely ergodic).  (The limit in
(\ref {L}) exists by subadditivity, which also
shows that $\lim$ may be replaced by $\inf$.)

\begin{rem} \label {sub}

By subadditivity,
$\frac {1} {2^k} \int \ln \|A^{(E)}_{2^k}(x)\| d\mu(x)$ is a decreasing
sequence converging to $L(E)$.  Allowing $E$ to take values in $\C$, we
conclude that $E \mapsto L(E)$ is the real part of a subharmonic function.

\end{rem}

\begin{lemma} \label {non-zero}

If $n \geq 2$, for every non-zero vector $z \in \R^2$, the derivative (with
respect to $E$) of the argument of $A^{(E,f,v)}_n(x) z$ is strictly negative.

\end{lemma}

\begin{pf}

Let $\rho_n(E,x,z)$ be the derivative (with respect to $E$)
of the argument of $A^{(E,f,v)}_n(x) z$.  It is easy to see that
$\rho_1(E,x,z)$ is strictly negative whenever $z$ is
not vertical, and it is zero if $z$ is vertical.  By the chain rule, for $n
\geq 2$, $\rho_n(E,x,z)=\sum_{i=1}^n \kappa_i
\rho_1(E,f^{i-1}(x),A^{(E,f,v)}_{i-1}(x) z)$, where
$\kappa_i$ are strictly positive (since $A^{(E,f,v)}_{n-i}(f^i(x))
\in \SL(2,\R)$ and hence preserves orientation).  Since either $z$ or
$A^{(E,f,v)}_1(x) z$ is non-vertical, the result follows.
\end{pf}

\subsubsection{} \label {l}

Let us endow the space $\HH$ of bounded
self-adjoint operators of $\ell^2(\Z)$
with the norm $\|\Phi\|=\sup_{\|u\|_2=1} \|\Phi(u)\|_2$, and the space of
compact subsets of $\R$ with the Caratheodory metric ($d(A,B)$ is the
infimum of all $r$ such that $A$ is contained in the $r$-neighborhood of $B$
and $B$ is contained in the $r$-neighborhood of $A$).  With respect to those
metrics, it is easy to see that the spectrum is a $1$-Lipschitz function of
$\Phi \in \HH$.  Since
the map $C^0(X,\R)$, $v \mapsto H_{f,v,x}$ is also $1$-Lipschitz, we
conclude that the spectrum of $H_{f,v,x}$ is a $1$-Lipschitz function of
$v \in C^0(X,\R)$.  It also follows that the spectrum of $H_{f,v,x}$ depends
continuously on $x$.

Since $H_{f,v,x}$ and $H_{f,v,f(x)}$ have obviously
the same spectrum, and $f$ is minimal, we conclude that the spectrum is
actually $x$-independent.  We will denote it $\Sigma(f,v)$.

\subsubsection{} \label {dense}

We say that $v$ is periodic (of period $n \geq 1$) if
$v(f^n(x))=v(x)$ for every $x \in X$.  If $v$ is a periodic potential, then
it is locally constant, hence for any compact subgroup $Y \subset X$
contained in a sufficiently small neighborhood of $\id$, the function $v$
is defined over $X/Y$.  If $v \in C^0(X,\R)$ and $Y \subset X$ is a
compact subgroup of
finite index, then we can define another potential $v^Y$ by convolution with
$Y$: $v^Y(x)=\int_Y v(y \cdot x) d\mu_Y$ where $\mu_Y$ is the Haar measure
on $Y$.  The potential $v^Y$ is then periodic.  Since there are
compact subgroups with finite index contained in arbitrarily small
neighborhoods of $\id$, this shows that
the set of periodic potentials is dense in $C^0(X,\R)$.

\subsubsection{} \label {per}

If $v$ is $n$-periodic then $\tr A^{(E,f,v)}_n(x)$ is
$x$-independent and denoted $\psi(E)$.  Then
$L(E,f,v)$ is the logarithm of the spectral radius of $A^{(E,f,v)}_n(x)$,
for any $x \in X$. 
This shows that the Lyapunov exponent is a continuous function of both the
potential and the energy {\it when one restricts considerations to
potentials of period $n$}.

\subsubsection{} \label {spectrum}

We will need some basic facts on the spectrum of periodic potentials, see
\cite {AMS}, \S 3, for a discussion with further references.

If $v$ is periodic of period $n$ the spectrum $\Sigma(f,v)$
of $H$ is the set of $E \in
\R$ such that $|\psi(E)| \leq 2$.  Thus for periodic potentials, we have
$\Sigma(f,v)=\{E \in \R,\, L(E,f,v)=0\}$.

The function $\psi$ is a polynomial of degree $n$.  It can be shown that
$\psi$ has $n$ distinct real roots and
its critical values do not belong to $(-2,2)$, moreover,
$E$ is a critical point of $\psi$ with $\psi(E)=\pm 2$ if and only if
$A^{(E,v,f)}_n(x)=\pm \id$.  From this one
derives a number of consequences about the structure of periodic spectra:
\begin{enumerate}
\item The set of all $E$ such that
$|\psi(E)|<2$ has $n$ connected components
whose closures are called bands,
\item If $E$ is in the boundary of some band, we obviously have
$\tr A^{(E,f,v)}_n(x)=\pm 2$,
\item Conversely, if $\tr A^{(E,f,v)}_n(x)=\pm 2$,
$E$ is in the boundary of some band, thus the spectrum is the
union of the bands,
\item If two different bands intersect then their common boundary point
satisfies $A^{(E,f,v)}_n(x)=\pm \id$.
\end{enumerate}

\subsubsection{}

We will need some simple estimates on the Lebesgue
measure of the bands and of the spectrum.

\begin{lemma} \label {meas}

Let $v$ be a periodic potential of period $n$.
\begin{enumerate}
\item The measure of each band is at most $\frac {2 \pi} {n}$,
\item Let $C \geq 1$ be such that
for every $E$ in the union of bands, there exists $x \in X$ and $k
\geq 1$ such that $\|A^{(E,f,v)}_k(x)\| \geq C$.  Then the total
measure of the spectrum is at most $\frac {4 \pi n} {C}$.
\end{enumerate}

\end{lemma}

\begin{pf}

If $E$ belongs to some band, $A^{(E,f,v)}_n(x)$ is conjugate in $\SL(2,\R)$
to a rotation: there exists $B^{(E)}(x) \in \SL(2,\R)$ such that $B^{(E)}(x)
A^{(E,f,v)}_n(x) B^{(E)}(x)^{-1} \in \SO(2,\R)$.  This matrix is
not unique, since $R B^{(E)}(x)$ has the same property for
$R \in \SO(2,\R)$, but this is the only ambiguity.  In particular, the
Hilbert-Schmidt norm squared $\|B^{(E)}(x)\|^2_\HS$
(the sum of the squares of the entries of the matrix of $B^{(E)}(x)$) is a
well defined function $b^{(E)}(x)$, which obviously satisfies
$b^{(E)}(f^n(x))=b^{(E)}(x)$.
This allows us to define an $x$-independent
function $\hat b(E)$ which is zero if $E$ does not belong to a band and for
$E$ in a band is given by
\be
\hat b(E)=\frac {1} {4 \pi n} \sum_{i=0}^{n-1} b^{(E)}(f^i(x)).
\ee
It turns out that $\hat b(E)$ is related to the {\it integrated density of
states} by the formula $N(E)=\int_{-\infty}^E \hat b(E) dE$.
As a consequence, we conclude that for any band
$I \subset \Sigma(f,v)$, $\int_I \hat b(E) dE=\frac {1} {n}$ (in particular
$\int_\R \hat b(E) dE=1$).  See \cite {AD2}, \S 2.4.1 for a discussion of
this point of view on the integrated density of states.

The first statement is then
an immediate consequence of $\hat b(E) \geq \frac {1}
{2 \pi}$ which in turn comes from the estimate $\|B\|^2_\HS \geq 2$,
$B \in \SL(2,\R)$.

For the second estimate, it is enough to show that for every $E$ in a band
we have $\hat b(E) \geq \frac {C} {4 \pi n}$.
Notice that
\begin{align}
B^{(E)}(f^k(x)) &A^{(E,f,v)}_k(x) A^{(E,f,v)}_n(x) A^{(E,f,v)}_k(x)^{-1}
B^{(E)}(f^k(x))^{-1}\\
\nonumber
&=B^{(E)}(f^k(x)) A^{(E,f,v)}_n(f^k(x))
B^{(E)}(f^k(x))^{-1}\in \SO(2,\R).
\end{align}
Thus
$B^{(E)}(f^k(x)) A^{(E,f,v)}_k(x)$ conjugates
$A^{(E,f,v)}_n(x)$ to a rotation so it coincides with $R B^{(E)}(x)$ for
some $R \in \SO(2,\R)$.  Thus
\be
C \leq \|A^{(E,f,v)}_k(x)\| \leq
\|B^{(E)}(f^k(x))^{-1}\| \|B(x)\|,
\ee
and there exists $y \in X$ (either $y=x$
or $y=f^k(x)$) such that $C \leq \|B^{(E)}(y)\|^2 \leq b^{(E)}(y)$.
It follows
that $\hat b(E) \geq \frac {C} {4 \pi n}$.
\end{pf}

\subsubsection{}

We conclude with a weak continuity result for the Lyapunov exponent.

\begin{lemma} \label {vn}

Let $v^{(n)} \in C^0(X,\R)$ be a sequence converging uniformly to $v \in
C^0(X,\R)$.  Then $L(E,f,v^{(n)}) \to L(E,f,v)$ in $L^1_\loc$.

\end{lemma}

\begin{pf}

This follows from the proof of Lemma 1 of \cite {AD1}.
Indeed for every compact interval
$I \subset \R$, there exists a continuous function $g:I \to \R$,
non-vanishing in $\inter I$, such that
\be
\lim_{n \to \infty}
\int_I \max \{L(E,f,v^{(n)})-L(E,f,v),0\} g(E) dE=0
\ee
and
\be
\lim_{n \to \infty}
\int_I \min \{L(E,f,v^{(n)})-L(E,f,v),0\} g(E) dE=0
\ee
(see the last two equations in page 396 of \cite {AD1}).  The result
follows.
\end{pf}


\section{Proof of Theorem \ref {coupl}}

Fix some Cantor group $X$, and let $f:X \to X$ be a minimal translation.
Then the homomorphism $\Z \to X$, $n \mapsto f^n(\id)$ is injective with
dense image.  For simplicity of notation, we identify the integers with its
image under this homomorphism.

For a given potential $w \in C^0(X,\R)$ and $n \geq 1$,
we write $L(E,w)=L(E,f,w)$ for the Lyapunov
exponent with energy $E$ corresponding to the potential $w$.


Since $X$ is Cantor, there exists
a decreasing sequence of Cantor subgroups $X_k \subset X$ with finite index
such that $\cap X_k=\{0\}$.
Let $P_k$ be the set of potentials which are defined on $X/X_k$.
Potentials in $P_k$ are $n_k$-periodic where $n_k$ is the index of $X_k$.
If $w \in C^0(X,\R)$ is a periodic potential, then it belongs to some $P_k$. 
Let $P=\cup P_k$ be the set of periodic potentials (which is a dense subset
of $C^0(X,\R)$, see \S \ref {dense}).

For $n \geq 1$, we write $A^{(E,w)}_n(x)=A^{(E,f,w)}_n(x)$
for the $n$-step transfer matrix associated with the potential $w$ at $x$.
We also let $A^{(E,w)}_n=A^{(E,w)}_n(0)$.  The
spectrum will be denoted by $\Sigma(w)=\Sigma(f,w)$.

We will actually work with finite families $W$ of periodic potentials.
Here we allow for multiplicity of elements, so the number of elements in
$W$, denoted by $\# W$, may be larger than the number of distinct elements
of $W$.  For simplicity of notation, we will often treat $W$ as a set
(writing for instance $W \subset P$).
We write $L(E,W)=\frac {1} {\# W} \sum_{w \in W} L(E,w)$.  (More formally,
and generally, one could work with probability measures with compact support
contained in $P_k$ for all $k$ sufficiently large.)

%

The core of the construction is contained in the following two lemmas.

\begin{lemma} \label {start}

Let $B$ be an open ball in $C^0(X,\R)$, let $W \subset P \cap B$ be a finite
family of potentials, and let
$M \geq 1$.  Then there exists a sequence $W^n \subset P \cap B$ such that
\begin{enumerate}
\item $L(E,\lambda W^n)>0$ whenever $M^{-1} \leq |\lambda| \leq M$,
$E \in \R$,
\item $L(E,\lambda W^n) \to L(E,\lambda W)$ uniformly on compacts (as
functions of $(E,\lambda) \in \R^2$).
\end{enumerate}

\end{lemma}

\begin{lemma} \label {induction}

Let $B$ be an open ball in $C^0(X,\R)$ and let $W \subset P \cap B$ be a
finite family of potentials.
Then for every $K$ sufficiently large, there exists
$W_K \subset P_K \cap B$ such that
\begin{enumerate}
\item $L(E,\lambda W_K) \to L(E,\lambda W)$
uniformly on compacts (as functions of $(E,\lambda) \in \R^2$),
\item The diameter of $W_K$ is at most $n_K^{-10}$,
\item For every $\lambda \in \R$, if
$\inf_{E \in \R} L(E,\lambda W) \geq \delta \# W n_k$ then
for every $w \in W_K$, $\Sigma(\lambda w)$ has Lebesgue measure at most
$e^{-\delta n_K/2}$.
\end{enumerate}

\end{lemma}

Before proving the lemmas, let us conclude the proof of Theorem \ref
{coupl}.  First we combine both lemmas:

\begin{lemma} \label {joining}

Let $B \subset C^0(X,\R)$ be an open ball and let $W \subset P \cap B$ be a
finite family of potentials. 
Then for every $M \geq 1$, there exist $\delta>0$,
an open ball $B'$ with closure contained in $B$,
with diameter at most $M^{-1}$ and
$W' \subset P \cap B'$ such that
\begin{enumerate}
\item $|L(E,\lambda W')-L(E,\lambda W)|<M^{-1}$ for $|E|,|\lambda| \leq M$,
\item $L(E,\lambda W')>\delta$ for every
$M^{-1} \leq |\lambda| \leq M$ and $E \in \R$,
\item For every $w \in B'$ and $M^{-1} \leq |\lambda| \leq M$ the Lebesgue
measure of $\Sigma(\lambda w)$ is at most $M^{-1}$.
\end{enumerate}

\end{lemma}

\begin{pf}

First apply Lemma \ref {start} to find some $\tilde W
\subset P \cap B$ such that
$L(E,\lambda \tilde W)>0$ for
every $E \in \R$ and $M^{-1} \leq |\lambda| \leq M$ (it is easy to see that
\be \label {e4}
L(E,\lambda w) \geq 1 \text { if } |E| \geq \|\lambda w\|+4,
\ee
so this is
really a statement about bounded energies which follows from Lemma \ref
{start}), and
$|L(E,\lambda \tilde W)-L(E,W)|<M^{-1}/4$ for every $|E|,|\lambda| \leq M$. 
By continuity of the Lyapunov exponent for periodic potentials (\S \ref
{per}) and compactness (and (\ref {e4}) to take care of large energies),
we conclude that there exists $\delta>0$
such that $L(E,\lambda \tilde W)>2\delta$ for every $E \in \R$ and $M^{-1}
\leq |\lambda| \leq M$.

Let us now apply Lemma \ref {induction} to $W=\tilde W$
and let $W'=W_K$ for $K$ large.  Then $W'$ is contained in a ball $B'
\subset B$ with diameter $n_K^{-10}<M^{-1}$ centered around some $w' \in
W'$.  Both estimates on $L(E,\lambda W')$ are clear from the statement of
Lemma \ref {induction}
(using again (\ref {e4}) for large $|E|$).
To estimate the measure of $\Sigma(\lambda w)$ for
$w \in B'$, we notice that $\Sigma(\lambda w)$
is contained in a $M n_K^{-10}$ neighborhood of
$\Sigma(\lambda w')$ (by $1$-Lipschitz
continuity of the spectrum, see \S \ref {l}).
Using that $\Sigma(\lambda w')$
has at most $n_K$ connected components and has measure at most $e^{-\delta
(\# \tilde W n_k)^{-1} n_K/2}$, the result follows.
\end{pf}

Given an open ball $B_0 \subset C^0(X,\R)$ and $W_0 \subset P \cap B_0$,
and $\epsilon_1>0$, we can proceed by induction,
applying the previous lemma, to define, for every $i \geq 1$,
open balls $B_i$ with $\overline B_i \subset B_{i-1}$, finite families of
periodic potentials $W_i \subset P \cap B_i$, and constants
$0<\delta_i<1$ and
$\epsilon_{i+1}=\min \{\epsilon_i,\delta_i\}/10$ such that
\begin{enumerate}
\item $L(E,\lambda W_i) \geq \delta_i$ for
$E \in \R$ and $\epsilon_i \leq |\lambda| \leq \epsilon^{-1}_i$,
\item $|L(E,\lambda W_i)-L(E,\lambda W_{i-1})|<\epsilon_i$
for $|E|,|\lambda| \leq \epsilon_i^{-1}$,
\item for every $w \in B_i$ and $\epsilon_i \leq |\lambda| \leq
\epsilon^{-1}_i$, $\Sigma(\lambda w)$ has measure at most $\epsilon_i$,
\end{enumerate}
Then the common element $w_\infty$
of all the $B_i$ is such that $\Sigma(\lambda w_\infty)$ has zero Lebesgue
measure for every $\lambda \neq 0$.  Notice that $L(E,\lambda W_i)$
converges uniformly on compacts to a continuous function, positive
if $\lambda \neq 0$, which by general considerations must coincide with
$L(E,\lambda w_\infty)$.
Indeed, if $w_n \to w$
then $L(E,w_n) \to L(E,w)$ in $L^1_\loc$ by Lemma \ref {vn}.
So $L(E,\lambda w_\infty)$ coincides almost
everywhere with $\lim L(E,\lambda W_i)$.
Since $E \mapsto L(E,\lambda w_\infty)$ is the real part of a subharmonic
function (see Remark \ref {sub}) and $E \mapsto \lim L(E,\lambda W_i)$
is continuous, they coincide everywhere.

Since $B_0$ was arbitrary, the
denseness claim of Theorem \ref {coupl} follows.

\subsection{Proof of Lemma \ref {start}}

Let $k$ be such that $W \subset P_k$.  For every $K>k$, choose $N_1(K)>0$
such that if $|E| \leq K$, $|\lambda| \leq K$, $w \in W$ and
$w' \in P_K$ are such that $w'$ is $\frac {2
n_k+1} {N_1(K)}$ close to $w$ then $|L(E,\lambda w')-L(E,\lambda w)|<\frac
{1} {K}$.  Here we use the continuity of the Lyapunov exponent for periodic
potentials, see \S \ref {per}.

For $w \in W$, $K>k$, $1 \leq j \leq 2 n_k+1$,
we define potentials $w^{K,j} \in P_K$ by
\be
w^{K,j}(i)=w(i),\, 0 \leq i \leq n_K-2 \text { and }
w^{K,j}(n_K-1)=w(n_K-1)+\frac {j} {N_1(K)}.
\ee
(This uniquely defines $w^{K,j}$ by periodicity.)

\begin{claim}

For every $\lambda \neq 0$, $K>k$ there exists $1 \leq j \leq 2
n_k+1$ such that $\Sigma(\lambda w^{K,j})$ has exactly $n_K$
components.

\end{claim}

\begin{pf}

Recall that for every $w' \in P_m$, there
exist exactly $2 n_m$ values of $E$ such that $\tr A^{(E,w')}_{n_m}=\pm
2$, if one counts the {\it exceptional}
energies such that $A^{(E,w')}_{n_m}=\pm \id$ with
multiplicity $2$, see \S \ref {spectrum}.

For each $j$ such that $\Sigma(\lambda w^{K,j})$ does not have exactly
$n_K$ components, there exists at least one energy $E_j \in
\Sigma(\lambda w^{K,j})$ with
$A^{(E_j,\lambda w^{K,j})}_{n_K}=
\pm \id$.  Then
\be
A^{(E_j,\lambda w)}_{n_K}=
\pm \left (\bm 1 & \frac {-\lambda j} {N_1(K)} \\ 0 & 1 \em \right ).
\ee
But since $w$ is $n_k$-periodic, this means that
\be
A^{(E_j,\lambda w)}_{n_k}=
\pm \left (\bm 1 &
\frac {-\lambda j n_k} {N_1(K) n_K} \\ 0 & 1 \em \right ).
\ee
This implies, in particular, that $A^{(E_j,\lambda w)}_{n_k} \neq
A^{(E_{j'},\lambda w)}_{n_k}$ for $j \neq j'$, thus we
must also have $E_j \neq E_{j'}$ for $j \neq j'$.  But there can be
at most $2 n_k$ values of $E$ such that
$\tr A^{(E,\lambda w)}_{n_k}=\pm 2$.  Thus there must be some $1 \leq j
\leq 2 n_k+1$ such
that $\Sigma(\lambda w^{K,j})$ has exactly $n_K$ connected components.
\end{pf}

By the previous claim and compactness, there exists $\delta=\delta(W,K,M)>0$
such that for $w \in W$ and $M^{-1} \leq |\lambda| \leq M$,
there exists $1 \leq j=j(K,\lambda,w) \leq 2 n_k+1$
such that $\Sigma(\lambda w^{K,j})$ has $n_K$ components and the
measure of the smallest gap is at least $\delta$.
Choose an integer $N_2(K)$ with $N_2(K)>\frac {4 \pi M} {\delta n_K}$.

For $0 \leq l \leq N_2(K)$ and $w^{K,j}$ as above, let
$w^{K,j,l} \in P_K$ be given by
$w^{K,j,l}=w^{K,j}+\frac {4 \pi M l} {n_K N_2(K)}$.

\begin{claim}

For every $M^{-1} \leq |\lambda| \leq M$, $w \in W$, $K>k$,
\be
\cap_{0 \leq l \leq N_2(K)} \Sigma(\lambda
w^{K,j(K,\lambda,w),l})=\emptyset.
\ee

\end{claim}

\begin{pf}

Each of the connected
components of $\Sigma(\lambda w^{K,j})$ has measure at most
$\frac {2 \pi} {n_K}$, see Lemma \ref {meas}.
Since $N_2(K)>\frac {4 \pi M} {\delta n_K}$, for every
$E$ there exists at least some $l$ with $0 \leq l \leq N_2(K)$
such that $E-\lambda \frac {4 \pi M l} {n_K N_2(K)} \notin \Sigma(\lambda
w^{K,j})$, that is, $E \notin \Sigma(\lambda w^{K,j,l})$, which gives the
result.
\end{pf}

Let $W^K$ be the family obtained by collecting the
$w^{K,j,l}$ for
different $w \in W$, $1 \leq j \leq 2 n_k+1$ and $0 \leq l \leq N_2(K)$.
By the second claim, $L(E,\lambda W^K)>0$ for every $M^{-1}
\leq |\lambda| \leq M$ and $E \in \R$ (since $L(E,\lambda w)>0$ if $E \notin
\Sigma(\lambda w)$, see \S \ref {spectrum}).

To conclude, it is enough to show that
\be
\max_{1 \leq j \leq 2 n_k+1}
\max_{0 \leq l \leq N_2(K)} |L(E,\lambda w^{K,j,l})-L(E,\lambda w)| \to 0
\ee
uniformly on compacts of $(E,\lambda) \in \R^2$.
Write
\begin{align}
|L(E,\lambda w^{K,j,l})-L(E,\lambda
w)| \leq &|L(E,\lambda w^{K,j,l})-L(E-\lambda \frac {4 \pi M l} {n_K N_2(K)},
\lambda w)|\\
\nonumber
&+|L(E-\lambda \frac {4 \pi M l} {n_K N_2(K)},\lambda w)-L(E,\lambda w)|.
\end{align}
Then the first term in the right hand side is smaller than
$K^{-1}$ provided $K \geq |E|+4 \pi M^2$ (by the choice of $N_1(K)$), while
the second term in the right hand side is bounded by $\max_{w \in W}
\sup_{|t| \leq \frac {4 \pi M^2} {n_K}} |L(E+t,\lambda w)-L(E,\lambda w)|$
which converges to zero uniformly on compacts of $(E,\lambda) \in \R^2$ as
$K \to \infty$ (by continuity of the Lyapunov exponent for periodic
potentials, \S \ref {per}).

\subsection{Proof of Lemma \ref {induction}}

Assume that $W \subset P_k$, $n_k \geq 2$,
and let $K>k$ be large.  Order the elements
$w^1,...,w^m$ of $W$.  Let $r=[n_K/m n_k]$.

First consider a
potential $w \in P_K$ obtained as follows.  It is enough to define $w(l)$
for $0 \leq l \leq n_K-1$.  Let
$I_j=[j n_k,(j+1) n_k-1] \subset \Z$ and let
$0=j_0<j_1<...<j_{m-1}<j_m=n_K/n_k$ be a
sequence such that $j_{i+1}-j_i-r \in \{0,1\}$.  Given $0 \leq l \leq
n_K-1$, let $j$ be such that $l \in I_j$, let $i$ be such that
$j_{i-1} \leq j<j_i$ and let $w(l)=w^i(l)$.

For any sequence $t=(t_1,...,t_m)$
with $t_i \in \{0,...,r-1\}$, let $w^t \in P_K$ be the potential defined as
follows.  Let $0 \leq l \leq n_K-1$, and let $j$ be such that $l \in I_j$. 
If $j=j_i-1$ for some $1 \leq i \leq m$, we let
$w^t(l)=w(l)+r^{-20} t_i$.  Otherwise we let $w^t(l)=w(l)$.

Let $W_K$
be the family consisting of all the $w^t$.  The claimed diameter estimate is
obvious for large $K$.

Let us show that $L(E,\lambda W_K) \to L(E,\lambda W)$
uniformly on compacts.
It is enough to restrict ourselves
to compact subsets of $(E,\lambda) \in \R \times
(\R \setminus \{0\})$, since it is easy to see that
$L(E,\lambda w)-L(E,0) \to 0$ uniformly as $\|\lambda w\| \to 0$.

For fixed $E$ and $\lambda$, we write
\be
A^{(E,\lambda w^t)}_{n_K}=C^{(t_m,m)} B^{(m)} \cdots
C^{(t_1,1)} B^{(1)},
\ee
where $C^{(t_i,i)}=A^{(E-\lambda r^{-20} t_i,\lambda
w^i)}_{n_k}$ and $B^{(i)}=(A^{(E,\lambda
w^i)}_{n_k})^{j_i-j_{i-1}-1}$.  Notice that, for $E$ and $\lambda$ in a
compact set, the norm of the $C^{(t_i,i)}$-type matrices stays bounded as
$r$ grows, while the $B^{(i)}$ matrices may get large.

Find some cutoff $(\ln \ln r)^{-m} \leq c \leq (\ln \ln \ln r)^m/(\ln \ln
r)^m$ such that if $\|B^{(i)}\|<e^{c r}$ then
$\|B^{(i)}\|<e^{(\ln \ln \ln r)^{-1} c r}<
e^{(\ln \ln \ln r)^{-1} c n_K}$.  To see that this is possible,
notice that the union of the $m$ intervals
$(\ln \ln \|B^{(i)}\|-\ln r,\ln \ln
\|B^{(i)}\|-\ln r+\ln \ln \ln \ln r]$, $1 \leq i \leq m$,
must ommit at least one point in
$[-m \ln \ln \ln r,-m \ln \ln \ln r+m\ln \ln \ln \ln r]$, which can
be taken as $\ln c$.


Call $i$ good if $\|B^{(i)}\| \geq e^{c r}$.  If no $B^{(i)}$ is good, then
$L(E,\lambda W) \leq c \frac {r} {r-1}$ and
$L(E,\lambda W_K) \leq c \frac {rm} {n_K}+O(1/r)$.  In particular
$L(E,\lambda W_K)$ and $L(E,\lambda W)$ are
close, since $c=o(1)$ with respect to $r$.

So we can assume that there exists at
least one good $B^{(i)}$.  Let $i_1<...<i_d$ be the list of all good $i$.
Write
$A^{(E,\lambda w^t)}(0)=\hat C^{(d)} \hat B^{(d)}
\cdots \hat C^{(1)} \hat B^{(1)}$, where for $1 \leq j \leq d$ we let
$\hat C^{(j)}=C^{(t_{i_j},i_j)}$
and $\hat B^{(j)}=B^{(i_j)} D^{(j)}$, where we denote
$D^{(j)}=C^{(i_j-1,t_{i_j-1})} B^{(i_j-1)}
\cdots C^{(i_{j-1}+1,t_{i_{j-1}+1})}
B^{(i_{j-1}+1)}$ (denoting also $i_0=0$).

By the choice of the cutoff, we have $\|D^{(j)}\| \leq e^{cr/2}$ for $r$
large (uniformly on compacts of $(E,\lambda) \in \R^2$), so
$\|\hat B^{(j)}\| \geq
e^{c r/2}$.

\begin{claim}

As $r$ grows,
\be \label {prod}
\frac {1} {n_K} \sum_{j=1}^d \ln \|\hat B^{(j)}\| \to L(E,\lambda W)
\ee
uniformly on compacts of $E$ and $\lambda$.

\end{claim}

\begin{pf}

Notice that this is equivalent to showing that
\be
\frac {1} {n_K} \sum_{i=1}^m \ln \|B^{(i)}\| \to L(E,\lambda W)
\ee
(uniformly), which in turn is equivalent to
\be
\frac {1} {m} \sum_{i=1}^m \frac {1} {n_k(j_i-j_{i-1}-1)}
\ln \|B^{(i)}\| \to L(E,\lambda W)
\ee
(uniformly).
Thus it is enough to show that
\be \label {growth}
\frac {1} {n_k (j_i-j_{i-1}-1)}
\ln \|B^{(i)}\| \to L(E,\lambda w^i)
\ee
(uniformly).  But $B^{(i)}$ is
just the $j_i-j_{i-1}-1$ iterate of the matrix $A^{(E,\lambda w^i)}_{n_k}$,
whose spectral radius is precisely the exponential of
$n_k L(E,\lambda w^i)$.  But it is easy to see that $\|T^n\|^{1/n}$
converges to the spectral radius of $T$ uniformly on compacts of $T \in
\SL(2,\R)$.  This gives (\ref {growth}) and the result.
\end{pf}

For every $t$, we have the obvious upper bound
\be \label {w^t}
L(E,\lambda w^t) \leq \frac {1} {n_K}
\sum_{j=1}^d \ln \|\hat B^{(j)}\|+O(1/r),
\ee
and we will now be concerned with bounding $L(E,\lambda w^t)$ from below,
not for all $t$, but for a majority of them.

Let $s_j$ be the most
contracted direction of $\hat B^{(j)}$ and let $u_j$ be the image under
$\hat B^{(j)}$ of the most expanded direction.

Let us say that $t$ is $j$-nice, $1 \leq j \leq d$,
if the absolute value of the
angle between $\hat C^{(j)} u_j$ and $s_{j+1}$ is
at least $r^{-70}$ (with the convention that $j+1=1$
for $j=d$).

\begin{claim}

Let $r$ be sufficiently large, and let $t$ be $j$-nice.  If $z$ is a
non-zero vector making an angle at least $r^{-80}$ with $s_j$,
then $z'=\hat C^{(j)} \hat B^{(j)} z$ makes an angle at least
$r^{-80}$ with $s_{j+1}$ and $\|z'\| \geq \|\hat B^{(j)}\| r^{-100} \|z\|$.

\end{claim}

\begin{pf}

Let $0 \leq \theta \leq \pi/2$
be the angle between $z$ and $s_j$, and let $0 \leq \theta' \leq \pi/2$
be the angle between $z''=\hat B^{(j)} z$ and $u_j$.

The orthogonal projection of
$z''$ on $u_j$ has norm $\|z\| \|\hat B^{(j)}\| \sin \theta$. 
Since $\|\hat C^{(j)}\|$ stays bounded as $r$ grows, we conclude that
$\|z'\| \geq \|\hat B^{(j)}\| r^{-100} \|z\|$.

On the other hand, $\tan \theta' \tan \theta=\|\hat B^{(j)}\|^{-2}$.  Since
$\|\hat B^{(j)}\| \geq e^{cr/2} \geq r^{400}$ for $r$ large, it follows that
$\theta'<r^{-100}$.  The boundedness of $\hat C^{(j)}$ again implies that
the angle between $z'$ and $\hat C^{(j)} u_j$ is at most $r^{-90}$.  Since
$t$ is $j$-nice, $z'$ makes an angle at least $r^{-80}$ with $s_{j+1}$.
\end{pf}

It follows that if $t$ is very nice in the sense that it is
$j$-nice for every $1 \leq j \leq d$, then if $z$
is a non-zero vector making an angle at least $r^{-80}$ with $s_1$ then
$z'=A^{(E,\lambda w^t)}_{n_K} z$ also makes an angle at least $r^{-80}$ with
$s_1$, and moreover $\|z'\|/\|z\| \geq
\prod_{j=1}^d r^{-100} \|\hat B^{(j)}\|$.
By (\ref {prod}) and (\ref {w^t}), it follows that
$L(E,\lambda w^t)-L(E,\lambda W) \to 0$ as $r$ grows, at least for
very nice $t$.

To conclude the estimate on the Lyapunov exponent, it is thus enough to show
that most $t$ are nice, in the sense that for every $\epsilon>0$,
for every $r$ sufficiently large,
the set of $t \in \{0,...,r-1\}^m$ which are
not very nice has at most $\epsilon r^m$ elements.  A more precise estimate
is provided below.

\begin{claim}

For every $r$ sufficiently large, the set of $t$ which are not very nice has
at most $m r^{m-1}$ elements.

\end{claim}

\begin{pf}

We will show in fact that, for every $1 \leq j \leq d$, if for every
$1 \leq k \leq m$ with $k \neq i_j$
one chooses $t_k \in \{0,...,r-1\}$, there exists at
most one ``exceptional'' $t_{i_j} \in \{0,...,r-1\}$ such that
$t=(t_1,...,t_m)$ is not $j$-nice.  Thus the set of $t$ which are not
$j$-nice has at most $r^{m-1}$ elements and the estimate follows.

Once $t_k$ is fixed for $1 \leq k \leq m$ with $k \neq i_j$, both $u_j$ and
$s_{j+1}$ become determined, but
$\hat C^{(j)}=C^{(t_{i_j},i_j)}=A^{(E-\lambda r^{-20} t_{i_j},
\lambda w^{i_j})}_{n_k}$
depends on $t_{i_j}$.

Since $n_k \geq 2$, we can apply Lemma \ref {non-zero} to conclude that
for any non-zero vector $z \in \R^2$,
the derivative of the argument of the vector $A^{(E',\lambda
w^{i_j})}_{n_k} z$ as a function of $E'$ is strictly negative, and hence
bounded away from zero and infinity, uniformly on $z$ and on
compacts of $(E',\lambda) \in \R^2$, and
independently of $r$.

If $r$ is sufficiently large, we conclude that
for every $0 \leq l \leq r-2$, there exists
a rotation $R_l$ of angle $\theta$ with $r^{-21}<\theta<r^{-19}$ such that
$C^{(l+1,i_j)} u_j=R_l C^{(l,i_j)} u_j$.  It immediately follows that
there exists at most one choice of $0 \leq t_{i_j} \leq r-1$ such that
$C^{(t_{i_j},i_j)} u_j$ has angle at most $r^{-90}$ with $s_{j+1}$,
as desired.
\end{pf}  

\comm{
We conclude that, fixing all cordinates of $t$ except the
ones of the form $t_{i_j}$ (this data is enough to define all the
$\hat B^{(j)}$, and hence all the $u_j$ and $s_j$),
there exists $\hat t_j \in \{0,...,r-1\}$ such that if
$t_{i_j} \neq \hat t_j$ then
$\hat C^{(j)} \hat B^{(j)} \cdot u_j$ has angle at least $r^{-100}$ with
$s_{j+1}$, with the convention that $s_{d+1}=s_1$.
Thus, for the majority of the choices of $t$
we have such an angle estimate, which yields, by
a straightforward computation
$|L(E,\lambda w^t)-L(E,\lambda w)| \to 0$ as $r \to \infty$, uniformly on
compacts of $E$ and $|\lambda| \geq \epsilon$.
}

We now estimate the measure of the spectrum.
\comm{
Clearly for every $w \in W_K$,
the integrated density of states $N$
has an analytic density in the interior of
the spectrum.  Thus it is enough to show that
$\frac {d N} {dE} \geq e^{\delta n_K/2}$
in the interior of the spectrum.  For fixed $E$ we have
\be
\frac {dN} {dE}=\frac {1} {4 \pi n_K} \sum_{j=0}^{n_K-1} \|B(j)\|^2_\HS,
\ee
where $\| \cdot \|_\HS$ stands for the Hilbert-Schmidt norm and the $B(j)$
are any matrices satisfying $B(j+1) A^{(E,\lambda
w)}_1(j) B(j)^{-1} \in \SO(2,\R)$
(which always exist in the interior of the spectrum).  This is a
reformulation of the usual formulas for periodic potentials,
see for instance \S 2.4.1 of \cite {AD2} for more details.
From the definition of $B(j)$, we have that for
$0 \leq a,b \leq n_K-1$ with $a<b$,
$\|B(a)\| \|B(b)\| \geq
\|A^{(E,\lambda w)}_{b-a}(a)\|$.
}
Let $w^i \in W$ be such that $L(E,\lambda
w^i) \geq \delta n_k m$.  Then
$\|A^{(E,\lambda w^t)}_{(r-1) n_k}((j_{i-1}) n_k)\| \geq e^{\delta m (r-1)
n_k^2}$.  Since $E$ is arbitrary, we can apply Lemma \ref {meas} to conclude
that the measure of the spectrum is at most
$4 \pi n_K e^{-\delta m (r-1) n_k^2} \leq e^{-\delta n_K/2}$ for $r$ large.
\comm{
We conclude that
\be
\frac {1} {4 \pi n_K} \sum_{j=0}^{n_K-1} \|B(j)\|^2_\HS
\geq e^{\delta n_K/2}
\ee
for $r$ large.
}
The result follows.

\end{document}